\documentclass[11pt]{article}
\usepackage{amsmath,latexsym,mathrsfs}
\usepackage{amssymb,amsfonts,amscd,amssymb,amsthm}
\usepackage{graphicx,graphics,verbatim,url}
\usepackage{pstricks}
\usepackage{color}

\numberwithin{equation}{section}
\newtheorem{theorem}{Theorem}[section]
\newtheorem{lemma}[theorem]{Lemma}
\newtheorem{corollary}[theorem]{Corollary}
\newtheorem{proposition}[theorem]{Proposition}

\newcommand{\R}{{\mathbb R}}

\newcommand{\mP}{\mathbb P}

\newcommand{\inpro}[2]{\left\langle{#1},{#2}\right\rangle}

\newcommand{\norm}[2]{\|{#1}\|_{#2}}

\newcommand{\snorm}[2]{\left|{#1}\right|_{#2}}

\newcommand{\G}{\Gamma}

\newcommand{\vecx}{\boldsymbol{x}}
\newcommand{\vecy}{\boldsymbol{y}}

\DeclareMathOperator{\dom}{{Dom}}

\DeclareMathOperator{\spann}{{span }}

\newcommand{\cD}{{\cal D}}

\newcommand{\cN}{{\cal N}}

\newcommand{\cT}{{\cal T}}

\newcommand{\gotoo}{\longrightarrow}

\newcommand{\wth}{\widehat}

\newcommand{\ol}{\overline}

\newcommand{\pa}{\partial}

\newcommand{\Proof}{{\bf Proof. }}

\newcommand{\ds}{\, ds}

\newcommand{\dx}{\, dx}

\newcommand{\nn}{\nonumber}

\date{}
\title{A mixed method for Dirichlet problems with radial basis
functions
\thanks{Partially supported by {CONICYT}-Chile through {FONDECYT}-project 1110324 and
Anillo ACT1118 (ANANUM), and UNSW FRG Grant number PS24436.}}
\author{Norbert Heuer
        \thanks{Facultad de Matem\'aticas,
                Pontificia Universidad Cat\'olica de Chile,
                Santiago, Chile.
                mailto:
                {\tt nheuer@mat.puc.cl}
               }
        \and
        Thanh Tran
        \thanks{School of Mathematics and Statistics,
                The University of New South Wales,
                Sydney 2052, Australia.
                mailto:
                {\tt Thanh.Tran@unsw.edu.au}
               }
}
\begin{document}
\maketitle
\begin{abstract}
We present a simple discretization by radial basis functions
for the Poisson equation with Dirichlet boundary condition.
A Lagrangian multiplier using piecewise polynomials
is used to accommodate the boundary condition. This
simplifies previous attempts to use radial basis functions
in the interior domain to approximate the solution 
and on the boundary to approximate the multiplier,
which technically requires that the mesh norm 
in the interior domain is
significantly smaller than that on the boundary.
Numerical experiments confirm theoretical results.


\end{abstract}

\section{Introduction}\label{sec:int}

Radial basis functions (RBFs) have been used successfully
\cite{Wen99} to solve partial differential equations with
Neumann or Robin boundary conditions. When Dirichlet
conditions are considered they must be approximated in an
appropriate way.
In the case of RBFs this causes analytical difficulties since
traces of such functions have representations that change
continuously with the position of their centers. Therefore,
corresponding discrete spaces (on the boundary) do not have a
clear structure that could be used for analysis. Such difficulties
do not appear when treating natural boundary conditions.

Lagrangian multipliers {consisting of RBFs}
have been studied in \cite{DuanTan2006}.
Due to a weaker inverse property of RBFs
as compared to that of finite elements,
the condition imposed on the mesh norms in the
interior domain and on the boundary is too restrictive. In
this paper we suggest to use RBFs in the domain and finite
elements on the boundary. We also {analyze the influence of
scaling of RBFs on the error estimate.} {Scaling of}
RBFs {can be used} to avoid too much overlap
{which is essential for the conditioning of stiffness matrices}.

The idea to couple RBFs with finite elements for the
Lagrangian multiplier is
proposed in~\cite{HeuerTran2012} where we solve a
hypersingular integral equation. In that paper, in order to
deal with the property that the solution of the integral
equation
can be extended by zero, we have to consider the Lagrangian
multiplier on an extended domain. In the current situation
it suffices to use multipliers on the boundary of the domain
where the problem is set. Moreover, not as
in~\cite{HeuerTran2012} where the Lagrangian multiplier does
not have any physical meaning, in the problem to be
considered in this paper it represents the normal derivative
of the solution. We prove
an inf-sup condition which allows an estimate for the
approximation of this multiplier; {we did not succeed in proving}
this condition for the problem considered in~\cite{HeuerTran2012}.

The paper is organized as follows. In Section~\ref{sec:mod
pro} we introduce the model problem and the mixed
variational formulation. In Section~\ref{sec:dis} the
finite-dimensional spaces and discretization are introduced.
Section~\ref{sec:app} is a revisit of approximation properties
of RBFs where we extend previous results so that they can be
used in the current study. Section~\ref{sec:opt} presents
the main result of the paper, namely {an a priori error
estimate} for the approximation.
This analysis is carried out after we prove
{discrete ellipticity of the Dirichlet bilinear form (in the
case that no mass term is present) and}
a discrete inf-sup condition {of the bilinear form}
involving the Lagrangian multiplier.
Section~\ref{sec:num} {presents} numerical experiments
{that confirm the error estimate}.

Throughout the paper, the notation $a \lesssim b$ indicates
that there exists a constant
$C>0$ independent of discretization or scaling parameters $h_X$, $k$, $r$
and involved functions (except where otherwise noted) such 
that $a \le C b$.
Similarly we use $a\gtrsim b$, and $a\simeq b$ means that
$a \lesssim b \lesssim a$.

\section{Model problem and mixed formulation}\label{sec:mod pro}

Consider the model problem
\begin{equation}\label{equ:mod pro}
\begin{aligned}
-\Delta u + \kappa u &= f \quad\text{ in }\Omega, \\
u &= g \quad\text{ on } \Gamma,
\end{aligned}
\end{equation}
where $\Omega\subset\R^d$ {($d=2,3$)} is a {polygonal ($d=2$) or 
polyhedral ($d=3$)}
domain with boundary $\Gamma$, $\kappa\ge0$ is a constant, and where
$f\in L_{2}(\Omega)$ and $g\in H^{1/2}(\Gamma)$
are given functions. 
{When $\kappa=0$ then $f$ and $g$ must satisfy the usual
compatibility condition.}

The solution $u$ of \eqref{equ:mod pro} will be found in
the weak sense by using a standard mixed variational
formulation.
Defining
\begin{equation}\label{equ:bil for}
\begin{aligned}
a(v,w)
&=
\int_{\Omega} (\nabla v \cdot \nabla w + \kappa vw) \dx
\quad\forall v,w\in H^1(\Omega), \\
b(v,\mu)
&=
\int_{\Gamma} v\mu \ds
\qquad\qquad\qquad\quad\forall v\in H^1(\Omega), \ \mu\in H^{-1/2}(\G),
\end{aligned}
\end{equation}
we can easily see that there exists a constant $C>0$ such that
\begin{equation}\label{equ:bdd}
\begin{aligned}
a(v,w)
&\le
C\norm{v}{H^1(\Omega)}\norm{w}{H^1(\Omega)}
\quad\forall v,w\in H^1(\Omega), \\
b(v,\mu)
&\le
C\norm{v}{H^1(\Omega)}\norm{\mu}{H^{-1/2}(\Gamma)}
\quad\forall v\in H^1(\Omega), \ \mu\in H^{-1/2}(\G).
\end{aligned}
\end{equation}
A variational formulation of \eqref{equ:mod pro} is
formulated as: Find
$(u,\lambda)\in H^1(\Omega)\times H^{-1/2}(\G)$
satisfying
\begin{equation}\label{equ:var for}
\begin{aligned}
a(u,v) + b(v,\lambda) &= F(v) \quad\forall v\in H^1(\Omega), \\
         b(u,\mu)     &= G(\mu) \quad\forall \mu\in
H^{-1/2}(\G),
\end{aligned}
\end{equation}
where
\[
F(v) = \displaystyle\int_{\Omega} fv\dx
\quad\text{and}\quad
G(\mu)=\displaystyle\int_{\Gamma}g\mu\ds.
\]

The following result is well known.
\begin{proposition}\label{pro:exi uni}
There exists a unique solution
$(u,\lambda)\in H^1(\Omega)\times H^{-1/2}(\G)$
to the problem \eqref{equ:var for}.
Moreover, $\lambda=\pa u/\pa n$ where $n$ is the outward
normal vector on $\G$.
\end{proposition}

\section{Discretization with RBFs and finite elements}\label{sec:dis}

We first define the finite-dimensional space that
approximates $u$ in~\eqref{equ:var for}.
Let $\Phi:\R^d\to\R$ be a radial basis function {whose Fourier
transform $\wth\Phi$ satisfies}
\begin{equation}\label{equ:Phi hat}
\wth\Phi(\omega)
\simeq
(1 + \snorm{\omega}{}^2)^{-\tau},
\quad
\omega\in\R^d,
\end{equation}
where $\tau>d/2$.
We consider the scaled radial basis functions
\[
\Phi_r(\vecx)
:=
r^{-d}\Phi(\vecx/r),
\quad r>0, \ \vecx\in\R^d,
\]
so that
\begin{equation}\label{equ:Phi r hat}
\wth\Phi_r(\omega)
\simeq
(1 + r^{2} \snorm{\omega}{}^2)^{-\tau},
\quad
\omega\in\R^d.
\end{equation}
The native space associated with $\Phi_r$ is defined by
\[
\cN_{\Phi_r}
=
\Bigl\{v\in\cD'(\R^d) \, : \,
\int_{\R^d} \frac{|\wth v(\omega)|^2}{\wth\Phi_r(\omega)}
\, d\omega < \infty \Bigr\}
\]
where $\cD'(\R^d)$ is the space of distributions defined in $\R^d$.
This space is equipped with an inner product and a norm defined by
\[
\inpro{v}{w}_{\Phi_r}
=
\int_{\R^d} \frac{\wth v(\omega)\overline{\wth w(\omega)}}{\wth\Phi_r(\omega)} \
d\omega
\quad\text{and}\quad
\norm{v}{\Phi_r} :=
\Bigl(\int_{\R^d} \frac{|\wth v(\omega)|^2}{\wth\Phi_r(\omega)}
\, d\omega\Bigr)^{1/2}.
\]
Under the assumption~\eqref{equ:Phi hat}, the native space
$\cN_{\Phi_r}$ is isomorphic to
the Sobolev space $H_r^{\tau}(\R^d)$ with equivalent norm
$\norm{\wth v(\cdot)(1+r^2\snorm{\cdot}{}^2)^{\tau/2}}{L_2(\R^d)}$.

Given a set of quasi-uniform centers
$X=\{\vecx_1,\ldots,\vecx_N\}\subset\Omega$ with mesh norm
$h_X = \sup_{\vecx\in\Omega}\min_{1\le j\le N}
\norm{\vecx-\vecx_j}{2}$,
we define
\begin{equation}\label{equ:HX}
H_{X,r}
=
\spann\{\Phi_1,\ldots,\Phi_N\},
\end{equation}
where $\Phi_i(\vecx) = \Phi_r(\vecx-\vecx_i)$ for $i=1,\ldots,N$.
(Note that since the nodes can be near to or even on the boundary
$\Gamma$ of $\Omega$, the supports of the scaled radial
basis functions are not necessarily subsets of~$\Omega$.)
The solution $u$ to~\eqref{equ:var for} will be approximated
by $u_X\in H_{X,r}$.

For the approximation of the Lagrangian multiplier~$\lambda$
we use functions (not necessarily continuous) which are
piecewise polynomials
of degree $p\ge0$ defined on
a {quasi-uniform} partition $\cT_k$ of the boundary $\G$ of $\Omega$:
\begin{equation}\label{equ:L1}
\Lambda_k :=
\{\mu:\G\to\R \ | \ \mu|_T \in \mP_p \ \forall
T\in\cT_k\}.
\end{equation}
Here, $k$ is the mesh size of $\cT_k$, and $\mP_p$ is the
space of polynomials of degree at most $p$.

Using these discrete spaces, the Galerkin scheme with radial
basis functions and Lagrangian multipliers for the
approximate solution of \eqref{equ:var for} is:
find $(u_X,\lambda_k) \in H_{X,r}\times\Lambda_k$ satisfying
\begin{equation}\label{equ:Gal equ}
\begin{aligned}
a(u_X,v) + b(v,\lambda_k) &= F(v) \quad\forall v\in H_{X,r}, \\
b(u_X,\mu) &= G(\mu) \quad\forall\mu\in\Lambda_k.
\end{aligned}
\end{equation}

\section{Approximation property of scaled RBFs}\label{sec:app}

For any integer $m\ge0$ and real $r>0$ we denote the norm of the
scaled Sobolev space $H_r^m(\Omega)$ by
\[
\norm{v}{H_r^m(\Omega)} :=
\left(\sum_{|\alpha|\le m} r^{2|\alpha|}
\norm{D^\alpha v}{L_2(\Omega)}^2\right)^{1/2}
\]
with multi-index $\alpha=(\alpha_1,{\ldots, \alpha_d})$ and
$|\alpha|=\alpha_1+{\cdots + \alpha_d}$.
We note that $H_r^0(\Omega)=H^0(\Omega)$.
For $s\in{(0,m)}$ {with integer $m>0$}  we define
the scaled interpolation space
\[
H_r^s(\Omega) := [H^0(\Omega),H_r^{m}(\Omega)]_\theta,
\]
where $\theta=s/m$.
{Here we employ the so-called real K-method; see
\cite{BerLof}.
The interpolation norm in $H_r^s(\Omega)$ can be represented as
follows. Let~$S_r:\dom(S_r)\subset H_r^m(\Omega)\to
H^0(\Omega)$ be an unbounded linear operator defined by
\[
\inpro{S_rv}{w}_{H^0(\Omega)}
=
\inpro{v}{w}_{H_r^m(\Omega)}
\quad\forall v \in \dom(S_r), \ \forall w \in H_r^m(\Omega).
\]
It is clear that $S_r$ is self-adjoint and positive. Thus there exists
$\Lambda_r := S_r^{1/2} : H_r^m(\Omega)\gotoo H^0(\Omega)$ satisfying
\[
\inpro{v}{w}_{H_r^m(\Omega)}
=
\inpro{\Lambda_rv}{\Lambda_rw}_{H^0(\Omega)}.
\]
The inner product and norm in $H_r^s(\Omega)$ can now be represented as
\begin{equation}\label{equ:Hrs}
\inpro{v}{w}_{H_r^s(\Omega)}
=
\inpro{\Lambda_r^{\theta}v}{\Lambda_r^{\theta}w}_{H^0(\Omega)}
\quad\text{and}\quad
\norm{v}{H_r^s(\Omega)}
=
\norm{\Lambda_r^\theta v}{H^0(\Omega)}.
\end{equation}
The Sobolev spaces $H_r^s(\Omega)$ form a Hilbert scale with the
following property.

\begin{lemma}\label{lem:int spa}
Let $s_1$ and $s_2$ be non-negative real numbers, and let
$s_0=(s_1+s_2)/2$. Then for any $v\in H_r^{s_1}(\Omega)\cap
H_r^{s_0}(\Omega)$ there holds
\[
\norm{v}{H_r^{s_1}{(\Omega)}}
=
\sup_{w\in\cD(\ol\Omega)\setminus\{0\}}
\frac{\inpro{v}{w}_{H_r^{s_0}(\Omega)}}{\norm{w}{H_r^{s_2}(\Omega)}}
\]
where $\cD(\ol\Omega)$ is the space of all functions which are
restrictions on $\ol\Omega$ of infinitely differentiable functions
in $\R^d$.
\end{lemma}
\Proof
Let $m$ be an integer not less than $\max\{s_1,s_2\}$.
{We may assume that $\max\{s_1,s_2\}>0$.} Then
\[
H_r^{s_i}(\Omega) = [H^0(\Omega),H_r^m(\Omega)]_{\theta_i},
\]
where $\theta_i=s_i/m$, $i=0,1,2$.
For any $v\in H_r^{s_1}(\Omega)\cap H_r^{s_0}(\Omega)$, there holds
\[
\norm{v}{H_r^{s_1}{(\Omega)}}
=
\sup_{z\in\cD(\ol\Omega)\setminus\{0\}}
\frac{\inpro{v}{z}_{H_r^{s_1}(\Omega)}}{\norm{z}{H_r^{s_1}(\Omega)}}.
\]
For each $z\in\cD(\ol\Omega)$ we define
$w=\Lambda_r^{\theta_1-\theta_2}z$.
{Then by noting~\eqref{equ:Hrs}, the relation
$\theta_0=(\theta_1+\theta_2)/2$, and the self-adjointness of
$\Lambda_r^\theta$, we obtain}
\[
\norm{z}{H_r^{s_1}(\Omega)}
=
\norm{\Lambda_r^{\theta_1}z}{H^{0}(\Omega)}
=
\norm{\Lambda_r^{\theta_2}w}{H^{0}(\Omega)}
=
\norm{w}{H_r^{s_2}(\Omega)}
\]
and
\begin{align*}
\inpro{v}{z}_{H_r^{s_1}(\Omega)}
&=
\inpro{\Lambda_r^{\theta_1}v}{\Lambda_r^{\theta_1}z}_{H^0(\Omega)}
=
\inpro{\Lambda_r^{\theta_1}v}{\Lambda_r^{\theta_2}w}_{H^0(\Omega)} \\
&=
\inpro{\Lambda_r^{\theta_0}v}{\Lambda_r^{\theta_0}w}_{H^0(\Omega)}
=
\inpro{v}{w}_{H_r^{s_0}(\Omega)}.
\end{align*}
Therefore,
\[
\norm{v}{H_r^{s_1}{(\Omega)}}
=
\sup_{w\in\cD(\ol\Omega)\setminus\{0\}}
\frac{\inpro{v}{w}_{H_r^{s_0}(\Omega)}}{\norm{w}{H_r^{s_2}(\Omega)}}
\]
and the lemma is proved.\qed

For any $v{\in H_r^\tau(\Omega)}$, let $I_Xv$ denote its interpolant in the space $H_{X,r}$,
i.e., $I_Xv\in H_{X,r}$ satisfies
\[
I_Xv(\vecx_j) = v(\vecx_j), \quad j = 1, \ldots, N.
\]
{For a domain in two dimensions ($d=2$),} the 
following approximation property is proved in
\cite[Lemma 5.3]{HeuerTran2012}; see also the proof of this lemma.
{The same arguments apply also to the case $d=3$.}

\begin{lemma}\label{lem:HX}
Let assumption \eqref{equ:Phi hat} be satisfied.
Then for any $v\in H_r^{\tau}(\Omega)$ there holds
\[
\norm{v-I_Xv}{H_r^s(\Omega)}
\le
C(s,\tau) \left(\frac{h_{X}}{r}\right)^{\tau-s}
\norm{v-I_Xv}{H_r^{\tau}(\Omega)}
\le
C(s,\tau) \left(\frac{h_{X}}{r}\right)^{\tau-s}
\norm{v}{H_r^{\tau}(\Omega)}
\]
for $0 \le s \le \lfloor\tau\rfloor$ and $0 < r \le r_0$
{with $r_0>0$ arbitrary but fixed}.
\end{lemma}

When $v$ is smoother, the error bound
can be {extended} by using the
technique developed in \cite{Sch99}, \cite{Wen05}, and modified in
\cite{TraLeGSloSte07a} for a sphere.

\begin{lemma}\label{lem:T}
Assume that \eqref{equ:Phi hat} holds. Let $T$ be an operator defined
by
\[
T\psi(\vecx) :=
\int_{\R^d} \Phi_r(\vecx-\vecy)\psi(\vecy) \ d\vecy,
\quad\vecx\in\R^d.
\]
Then for any $s\in\R$
\begin{enumerate}
\item
$T$ is an isomorphism from $H_r^{s-\tau}(\R^d)$ onto
$H_r^{s+\tau}(\R^d)$ and satisfies
\[
\norm{T\psi}{H_r^{s+\tau}(\R^d)}
\simeq
\norm{\psi}{H_r^{s-\tau}(\R^d)};
\]
\item
For any $\psi\in H_r^{-\tau}(\R^d)$ and $\xi\in\cN_{\Phi_r}$ there
holds
\[
\inpro{T\psi}{\xi}_{\Phi_r}
=
\inpro{\psi}{\xi}_{H_r^0(\R^d)},
\]
i.e., $T$ is the adjoint of the embedding operator of the native space
$\cN_{\Phi_r}$ into $H_r^0(\R^d)$.
\end{enumerate}
\end{lemma}
\Proof
The lemma follows from
\[
\widehat{T\psi}(\omega)
=
\widehat{\Phi}_r(\omega) \widehat{\psi}(\omega)
\simeq
(1+r^{2} |\omega|^2)^{-\tau} \widehat{\psi}(\omega),
\quad\omega\in\R^d.
\]
\qed

With the help of the above lemma, {we now extend} the error bound
in Lemma~\ref{lem:HX} when $v$ is smoother.

\begin{lemma}\label{lem:HX imp}
Let the assumption \eqref{equ:Phi hat} be satisfied. Then for any
$s,t\in\R$ satisfying $0\le s \le \tau \le t \le 2\tau$, if $v\in
H_r^t(\Omega)$ then the following estimate holds
\[
\norm{v-I_Xv}{H_r^s(\Omega)}
\le
C(s,t,\tau)
\left(\frac{h_X}{r}\right)^{t-s}
\norm{v}{H_r^t(\Omega)}
\]
{for $0<r\le r_0$ with $r_0>0$ arbitrary but fixed}.
\end{lemma}
\Proof
Consider first the case when $t=2\tau$.
Let
\[
E : H_r^\sigma(\Omega) \to H_r^\sigma(\R^d)
\]
be an $r$-uniformly bounded extension operator
for any $\sigma>0$; cf. \cite[Lemma 5.1]{HeuerTran2012}.
For any $v\in H_r^{2\tau}(\Omega)$,
since $I_XEv = I_Xv = EI_Xv$
on $\Omega$, one finds that $Ev-I_XEv$ is an extension of $v-I_Xv$.
Therefore, the property that $I_X$ is an orthogonal projection in
$\cN_{\Phi_r}$ yields
\begin{align*}
\norm{v-I_Xv}{H_r^\tau(\Omega)}^2
&\le
\norm{Ev-I_XEv}{H_r^\tau(\R^d)}^2
\simeq
\norm{Ev-I_XEv}{{\Phi_r}}^2 \\
&\le
\inpro{Ev-I_XEv}{Ev}_{{\Phi_r}}.
\end{align*}
It follows from Lemma~\ref{lem:T} that there exists $\psi\in
H_r^0(\R^d)$ such that
\[
T\psi=Ev,
\quad
\norm{\psi}{H_r^0(\R^d)} \simeq \norm{Ev}{H_r^{2\tau}(\R^d)}
\]
and
\[
\inpro{Ev-I_XEv}{Ev}_{{\Phi_r}}
=
\inpro{Ev-I_XEv}{\psi}_{H_r^0(\R^d)}.
\]
Therefore
\begin{align}\label{equ:tau 0}
\norm{v-I_Xv}{H_r^\tau(\Omega)}^2
&\lesssim
\norm{Ev-I_XEv}{H_r^0(\R^d)}
\norm{Ev}{H_r^{2\tau}(\R^d)} \nn \\
&\lesssim
\norm{v-I_Xv}{H_r^0(\Omega)}
\norm{v}{H_r^{2\tau}(\Omega)}.
\end{align}
By applying Lemma~\ref{lem:HX} with $s=0$ and using the
above inequality we obtain
\begin{align*}
\norm{v-I_Xv}{H_r^0(\Omega)}^2
&\le
C(\tau)
\left(\frac{h_X}{r}\right)^{2\tau}
\norm{v-I_Xv}{H_r^\tau(\Omega)}^2 \\
&\le
C(\tau)
\left(\frac{h_X}{r}\right)^{2\tau}
\norm{v-I_Xv}{H_r^0(\Omega)} \norm{v}{H_r^{2\tau}(\Omega)}.
\end{align*}
Thus the required estimate is proved for $s=0$ and $t=2\tau$.
It also follows from~\eqref{equ:tau 0} that the required estimate
holds for $s=\tau$ and $t=2\tau$. By using interpolation we deduce
the estimate for $0\le s\le \tau$ and $t=2\tau$.
Interpolation between
\[
\norm{v-I_Xv}{H_r^s(\Omega)}
{\lesssim}
\left(\frac{h_X}{r}\right)^{2\tau-s}
\norm{v}{H_r^{2\tau}(\Omega)}
\]
and
\[
\norm{v-I_Xv}{H_r^s(\Omega)}
{\lesssim}
\left(\frac{h_X}{r}\right)^{\tau-s}
\norm{v}{H_r^{\tau}(\Omega)}
\]
yields the required estimate for
$0\le s\le \tau \le t\le 2\tau$, proving the lemma.
\qed

To extend Lemma~\ref{lem:HX imp} to include less smooth functions,
i.e. $v\in H_r^t(\Omega)$ for $0\le t < \tau$,
we use Lemma~\ref{lem:int spa}.

\begin{lemma}\label{lem:app pro}
Assume that \eqref{equ:Phi hat} holds. For any $v\in
H_r^t(\Omega)$ with $0\le t \le 2 \tau$ there exists $z_X\in
H_{X,r}$ satisfying
\begin{equation}\label{equ:st}
\norm{v-z_X}{H_r^s(\Omega)}
\le
C(s,t,\tau)
\left(\frac{h_X}{r}\right)^{t-s}
\norm{v}{H_r^t(\Omega)}
\end{equation}
for $0\le s \le \min\{t, \tau\}$
{and $0<r\le r_0$ with $r_0>0$ arbitrary but fixed}.
\end{lemma}
\Proof
We only need to prove the lemma for $0\le t < \tau$.
Consider first the case when $\tau/2 \le t < \tau$.
Then $0 \le 2t-\tau < t$.
Let $P_t : H_r^t(\Omega)\gotoo H_{X,r}$ be the projection
defined by
\begin{equation}\label{equ:Pt0}
\inpro{P_tv}{z}_{H_r^t(\Omega)}
=
\inpro{v}{z}_{H_r^t(\Omega)}
\quad\forall z\in H_{X,r}.
\end{equation}
It can be seen that
\begin{equation}\label{equ:Pt}
\norm{P_tv-v}{H_r^t(\Omega)}
\le
\norm{v}{H_r^t(\Omega)}.
\end{equation}
If $s\in[0,2t-\tau]$ then
$\tau\le 2t-s < 2\tau$. Hence, for any $w\in H_r^{2t-s}$ it follows
from Lemma~\ref{lem:HX imp} that
\begin{equation}\label{equ:w Iw}
\norm{w-I_Xw}{H_r^t(\Omega)}
\lesssim
\left(\frac{h_X}{r}\right)^{t-s}\norm{w}{H_r^{2t-s}(\Omega)}.
\end{equation}
By using Lemma~\ref{lem:int spa} and
\eqref{equ:Pt0}--\eqref{equ:w Iw} we obtain
\begin{align*}
\norm{P_tv - v}{H_r^s(\Omega)}
&=
\sup_{w\in\cD(\ol\Omega)\setminus\{0\}}
\frac{\inpro{P_tv-v}{w}_{H_r^t(\Omega)}}{\norm{w}{H_r^{2t-s}(\Omega)}} \\
&=
\sup_{w\in\cD(\ol\Omega)\setminus\{0\}}
\frac{\inpro{P_tv-v}{w-I_Xw}_{H_r^t(\Omega)}}{\norm{w}{H_r^{2t-s}(\Omega)}} \\
&\lesssim
\left(\frac{h_X}{r}\right)^{t-s} \norm{v}{H_r^t(\Omega)}.
\end{align*}
In particular, there holds
\[
\norm{P_tv - v}{H_r^{2t-\tau}(\Omega)}
\lesssim
\left(\frac{h_X}{r}\right)^{\tau-t}
\norm{v}{H_r^t(\Omega)}.
\]
Hence, for $s\in[2t-\tau,t]$ by
noting~\eqref{equ:Pt} and using interpolation
we obtain the required estimate, and thus
prove~\eqref{equ:st}
for $0\le s \le t$ and $\tau/2\le t <\tau$.

By successively considering the case $\tau/4 \le t < \tau/2$, then
$\tau/8 \le t < \tau/4$, etc., and using the same argument,
we finish the proof of the lemma.
\qed

We are now able to derive the approximation property of
$H_{X,r}$ in non-scaled Sobolev norms.

\begin{lemma}\label{lem:app pro non}
Assume that \eqref{equ:Phi hat} holds. For any $v\in
H_r^t(\Omega)$ with $0\le t \le 2 \tau$ there exists $z_X\in
H_{X,r}$ satisfying
\[
\norm{v-z_X}{H^s(\Omega)}
\le
C(s,t,\tau)
\frac{h_X^{t-s}}{r^t}
\norm{v}{H^t(\Omega)}
\]
for $0\le s \le \min\{t, \tau\}$
{and $0<r\le r_0$ with $r_0>0$ arbitrary but fixed}.
\end{lemma}
\Proof
First we note that since {$r\in(0,r_0]$},
for any function $v$ and any positive
integer $m$ there hold
\[
\norm{v}{H^0(\Omega)} = \norm{v}{H_r^0(\Omega)}
\quad\text{and}\quad
\norm{v}{H_r^m(\Omega)}
{\lesssim}
\norm{v}{H^m(\Omega)}
\le
r^{-m} \norm{v}{H_r^m(\Omega)}.
\]
By interpolation it follows that {for $0 \le s \le m$}
\[
\norm{v}{H_r^s(\Omega)}
{\lesssim}
\norm{v}{H^s(\Omega)}
\le
r^{-s} \norm{v}{H_r^s(\Omega)}.
\]
The required result is then a consequence of the above
inequalities and Lemma~\ref{lem:app pro}.
\qed

\section{{Error estimate}}\label{sec:opt}

{In this section we prove our main result, Theorem~\ref{the:con}, which establishes
quasi-optimal convergence of our mixed method and convergence orders. Of course,
approximation properties depend on the regularity of solutions. To keep things simple,
we assume that, for smooth data, we have standard elliptic regularity limited by
the smoothness of $\Gamma$. More precisely, let $\delta$ be such that the
solution $u$ of \eqref{equ:mod pro} for any $\kappa\ge 0$ and any 
sufficiently smooth
data $f$, $g$, satisfies
\begin{equation} \label{reg}
   \delta\in (1,2]:\quad u\in H^\delta(\Omega).
\end{equation}
We restrict our considerations to regularity no more than $H^2(\Omega)$ since we
are interested in non-smooth problems and to simplify results when approximation
spaces use piecewise polynomials of higher degree. Of course, in two dimensions
$\Omega$ being a polygon, there holds
\[
   \delta =
   \begin{cases}
      2 \quad & \text{if $\Omega$ is convex}, \\
      2\pi/\omega \quad & \text{if $\Omega$ is non-convex},
   \end{cases}
\]
with $\omega\in(\pi,2\pi)$ being the angle of the largest
re-entrant corner of the boundary $\Gamma$ in case of non-convex $\Omega$.
The characterization of $\delta$ for a polyhedral domain is a bit more involved and
not given here.

We use standard Babu{\v s}ka-Brezzi theory to prove the main result. In order to
do so we now prove ellipticity of the bilinear form $a(\cdot,\cdot)$ on
\[
   V_{X,r} := \{v\in H_{X,r};\; b(v,\mu)=0\quad \forall \mu\in\Lambda_k\}
\]
when $\kappa=0$ and an inf-sup condition for the bilinear form $b(\cdot,\cdot)$.
}

\begin{lemma} \label{lem:ell}
For $k$ sufficiently small there holds
\[
   |v|_{H^1(\Omega)} \gtrsim \norm{v}{H^1(\Omega)}
   \qquad\forall v\in V_{X,r}.
\]
\end{lemma}
\Proof
The proof is standard
(cf.~\cite[Lemma~4.5]{GatHeaHeu09})
and is given for convenience of the reader.
Let $v\in V_{X,r}$ be given and decomposed as $v=v_0+{D}$ with
${D}=|\Omega|^{-1}\int_\Omega v\dx$ so that $\int_\Omega v_0\dx=0$.
Here, $|\Omega|$ denotes the measure of $\Omega$. There holds
\begin{equation} \label{dec}
   \norm{v}{H^1(\Omega)}^2 = \norm{v_0}{H^1(\Omega)}^2 + |\Omega| {D}^2.
\end{equation}
Now let $\mu\in L_2(\G)\setminus\{0\}$ and its $L_2(\G)$-projection
$\Pi\mu$ onto $\Lambda_k$ be given. By duality, a standard approximation result
and the trace theorem, we find that there holds
\begin{align*}
   |b({D},\mu)| &= |b(v-v_0,\mu)| = |b(v,\mu-\Pi\mu) - b(v_0,\mu)|\\
   &\lesssim
   \norm{v}{H^{1/2}(\G)} k^{1/2} \norm{\mu}{L_2(\G)}
   +
   \norm{v_0}{L_2(\G)}\norm{\mu}{L_2(\G)}\\
   &\lesssim
   \Bigl(k^{1/2} \norm{v}{H^{1}(\Omega)} + \norm{v_0}{H^1(\Omega)}\Bigr) \norm{\mu}{L_2(\G)},
\end{align*}
that is,
\begin{equation} \label{est_d}
   |{D}| \simeq \sup_{\mu\in L_2(\G)\setminus\{0\}} \frac{|b({D},\mu)|}{\norm{\mu}{L_2(\G)}}
   \lesssim
   k^{1/2} \norm{v}{H^{1}(\Omega)} + \norm{v_0}{H^1(\Omega)}.
\end{equation}
Now, using \eqref{dec}, \eqref{est_d} and Poincar\'e-Friedrichs' inequality,
we find that
\begin{align*}
   \norm{v}{H^1(\Omega)}^2
   \lesssim
   \norm{v_0}{H^1(\Omega)}^2 + k \norm{v}{H^1(\Omega)}^2
   \lesssim
   |v|_{H^1(\Omega)}^2 + k \norm{v}{H^1(\Omega)}^2.
\end{align*}
Selecting $k$ small enough finishes the proof.
\qed

\begin{lemma}\label{lem:inf sup}
Assume that \eqref{equ:Phi hat} 
and the elliptic regularity \eqref{reg} hold. Suppose that
the parameters $h_X$, $k$ and $r$ are chosen such that
\begin{equation}\label{equ:K}
   K_t(h_X,k,r):=\frac{h_X^{t-1}}{k^{t-1}r^{t}}
   \quad\text{is sufficiently small}
\end{equation}
with $t=\delta$ when $\delta<2$, or
with some $t\in [1,2)$ when $\delta=2$.
Then there exists a positive constant $\alpha$, 
independent of $h_X$, $k$ and $r$,
except for {their} relation via $K_t$,
such that there holds
\[
\alpha\norm{\mu_k}{H^{-1/2}(\Gamma)}
\le
\sup_{v_X\in H_{X,r}\setminus\{0\}}
\frac{b(v_X,\mu_k)}{\norm{v_X}{H^1(\Omega)}}
\qquad\forall \mu_k\in \Lambda_k.
\]
\end{lemma}

\Proof
For any $\mu_k\in\Lambda_k$, consider the problem
\begin{equation}\label{equ:wkw}
\begin{alignedat}{2}
-\Delta w + w &= 0 & \quad &\text{ in }\Omega, \\
\frac{\pa w}{\pa n} &= \mu_k & \quad &\text{ on } \Gamma.
\end{alignedat}
\end{equation}
Since $\mu_k\in H^{{s}}(\G)$ {for any $s<1/2$},
there exists a unique {variational solution} $w\in H^{{1}}(\Omega)$
{of \eqref{equ:wkw} with regularity estimate (limited by assumption \eqref{reg}
depending on the geometry)}
\begin{equation}\label{equ:w mu}
   \norm{w}{H^{\min\{\delta,3/2+s\}}(\Omega)}
   \lesssim
   \norm{\mu_k}{H^{\min\{\delta-3/2,s\}}(\G)}
   \qquad\forall s\in [-1/2,1/2);
\end{equation}
see e.g. \cite{Gri}.
Moreover, it is shown in \cite[Theorem 2.7]{Bab73} that
\begin{equation}\label{equ:w muk}
b(w,\mu_k) \simeq \norm{\mu_k}{H^{-1/2}(\G)}^2.
\end{equation}
On the other hand, since ${1}<\delta \le 2<2\tau$ we can invoke
Lemma~\ref{lem:app pro non} to obtain, for some
$z_X\in H_{X,r}$,
\[
\norm{w-z_X}{H^1(\Omega)}
\le
{c(t)} \frac{h_X^{{t-1}}}{r^{{t}}} \norm{w}{H^{{t}}(\Omega)}
   \qquad\text{with}\quad \left\{
   \begin{array}{ll}
      t=\delta & \text{if } \delta<2,\\
      \forall t\in [1,2) & \text{if } \delta=2.
   \end{array}\right.
\]
By using \eqref{equ:w mu}, the inverse property
\[
\norm{\mu_k}{H^{{t-3/2}}(\G)}
\lesssim
{k^{1-t}}
\norm{\mu_k}{H^{-1/2}(\G)}
   {\qquad \forall t\in [1,2)},
\]
and the assumption~\eqref{equ:K},
we deduce {that}
\begin{equation}\label{equ:w zX}
\norm{w-z_X}{H^1(\Omega)}
   \lesssim
   \frac{h_X^{t-1}}{r^t k^{t-1}}
   \norm{\mu_k}{H^{-1/2}(\G)}
=
{K_{t}(h_X,k,r)} \norm{\mu_k}{H^{-1/2}(\G)}
\end{equation}
with $t$ as given (if $\delta=2$) or chosen (if $\delta<2$) for \eqref{equ:K}.
This inequality, assumption \eqref{equ:K} and 
\eqref{equ:w mu} {with $s=-1/2$} give
\begin{equation}\label{equ:zX mu}
\norm{z_X}{H^1(\Omega)}
\lesssim
\norm{\mu_k}{H^{-1/2}(\G)}.
\end{equation}
On the other hand \eqref{equ:bdd},
\eqref{equ:w muk} and \eqref{equ:w zX} yield 
{for $K_t(h_X,k,r)$ small enough}
\begin{align*}
b(z_X,\mu_k)
&=
b(w,\mu_k)+  b(z_X-w,\mu_k) \\
&\gtrsim
\norm{\mu_k}{H^{-1/2}(\G)}^2
-
\norm{\mu_k}{H^{-1/2}(\G)} \norm{w-z_X}{H^1(\Omega)} \\
&\gtrsim
\left(1-{K_t(h_X,k,r)}\right) \norm{\mu_k}{H^{-1/2}(\G)}^2
\gtrsim
\norm{\mu_k}{H^{-1/2}(\G)}^2,
\end{align*}
i.e.,
\[
\norm{\mu_k}{H^{-1/2}(\G)}
\lesssim
\frac{b(z_X,\mu_k)}{\norm{\mu_k}{H^{-1/2}(\G)}}.
\]
This together with~\eqref{equ:zX mu} yields
\[
\norm{\mu_k}{H^{-1/2}(\G)}
\lesssim
\frac{b(z_X,\mu_k)}{\norm{z_X}{H^{1}(\Omega)}}
\le
\sup_{v_X\in H_{X,r}\setminus\{0\}}
\frac{b(v_X,\mu_k)}{\norm{v_X}{H^{1}(\Omega)}},
\]
proving the lemma.
\qed

{The following theorem is our main result. It proves the quasi-optimal
convergence of our mixed RBF approximation of the Dirichlet problem with
finite element Lagrangian multiplier.}

\begin{theorem}\label{the:con}
Let us assume that \eqref{equ:Phi hat} and \eqref{equ:K} hold and
consider radius parameters $r>0$ which are bounded. In the case that $\kappa=0$
in \eqref{equ:mod pro}, we additionally assume that $k$ is small enough. Then
there exists a unique solution $(u_X,\lambda_k)\in H_{X,r}\times\Lambda_k$
to the problem \eqref{equ:Gal equ}. Moreover, 
let $f$ and $g$ be sufficiently smooth so that
$(u,\lambda)$ is the solution to \eqref{equ:var for}
with $u\in H^\delta(\Omega)$ and $\lambda\in H^{\delta-3/2}(\Gamma)$
with $\delta\in (1,2]$, cf.~\eqref{reg}. Then
\begin{align*}
\norm{u-u_X}{H^1(\Omega)}
+
\norm{\lambda-\lambda_k}{H^{-1/2}(\Gamma)}
&\lesssim
\inf_{v\in H_{X,r}} \norm{u-v}{H^1(\Omega)}
+
\inf_{\mu_k\in \Lambda_k}
\norm{\lambda-\mu_k}{H^{-1/2}(\Gamma)} \\
&\lesssim
\frac{h_X^{\delta-1}}{r^\delta} \norm{u}{H^\delta(\Omega)}
+
k^{\delta-1}\norm{\lambda}{H^{\delta-3/2}(\G)},
\end{align*}
where the {implicitly appearing} constants depend on the constant $C$ in
\eqref{equ:bdd}, $\alpha$ in Lemma~\ref{lem:inf sup},
{and the ellipticity constant of bilinear form $a(\cdot,\cdot)$.}
\end{theorem}
\Proof
{For $\kappa>0$, the first bound of quasi-optimal convergence
follows from standard Babu\v{s}ka-Brezzi theory
(cf. \cite[Corollary~12.5.18]{BreSco02}) by making use of the continuity
and ellipticity of the bilinear form $a(\cdot,\cdot)$, and continuity and
continuous as well as discrete inf-sup condition of the bilinear form
$b(\cdot,\cdot)$, cf. Lemma~\ref{lem:inf sup}.

In the case that $\kappa=0$ we use Lemma~\ref{lem:ell} to obtain
$V_{X,r}$-ellipticity of $a(v,v)=\int_\Omega |\nabla v|^2 \dx$.
Then the result follows the same way.}

To show the second estimate for $\kappa\ge 0$ we use
Lemma~\ref{lem:app pro non} and the approximation property
of~$\Lambda_k$.
\qed

Depending on the regularity of the solution, the error estimate
above can be simplified.

\begin{corollary} \label{cor:con}
Let the assumptions of Theorem~\ref{the:con} be satisfied. If $\delta>3/2$ then
\begin{align*}
\norm{u-u_X}{H^1(\Omega)}
+
\norm{\lambda-\lambda_k}{H^{-1/2}(\Gamma)}
&\lesssim
\Bigl(\frac{h_X^{\delta-1}}{r^\delta} + k^{\delta-1}\Bigr)
\norm{u}{H^\delta(\Omega)}.
\end{align*}
\end{corollary}

\Proof
In the case $u\in H^\delta(\Omega)$ with $\delta>3/2$, the normal derivative
$\lambda=\partial u/\partial n$ can be defined in the standard way (normal
component trace of weak gradient) so that
\[
   \norm{\lambda}{H^{\delta-3/2}(\G)} \lesssim \norm{u}{H^\delta(\Omega)}.
\]
The assertion then is direct consequence of Theorem~\ref{the:con}.
\qed

\section{Numerical results}\label{sec:num}

We consider the model problem \eqref{equ:mod pro} with $\Omega=(0,1)\times (0,1)$,
$\kappa=0$ and $f$, $g$ such that $u(x_1,x_2)=x_1^2+x_2^2$, i.e. $u\in H^\delta(\Omega)$
with $\delta=2$.
The nodes of $X$ are distributed uniformly on $\bar\Omega$ including nodes on the boundary.

We use scaled radial basis functions with the radial basis functions defined in
\cite{Wen_EEI98}.
We consider the two cases $\tau=1.5$ and $\tau=2.5$ which
correspond to $C^0$ and $C^2$-functions, respectively. They are rotations of
univariate polynomials of degrees 2 and 5, respectively.
We have implemented the method
by numerical integration with an overkill of number of integration nodes.

With $\delta=2$ assumption \eqref{equ:K} requires that the ratio
$\frac{h_X^{1-\epsilon}}{k^{1-\epsilon}r^{2-\epsilon}}$
be small enough for some~$\epsilon>0$.
We simply choose $k\simeq h_{X,r}/r$ (more precisely an integer approximation to
$1/k$ smaller than or equal to $1$ since the length of the sides of $\Omega$ is $1$).
In this way, for fixed $r$, $K_t$ is fixed and \eqref{equ:K} is not guaranteed.
However, our numerical results do not show stability problems (that might be caused
by a violation of the inf-sup condition) in the range of unknowns under consideration.

For fixed $r$, the error estimate derived in
Corollary~\ref{cor:con} gives an upper bound
\[
   \norm{u-u_X}{H^1(\Omega)}
   +
   \norm{\lambda-\lambda_k}{H^{-1/2}(\Gamma)}
   \lesssim
   h_X + k.
\]
In the graphs below we plot the individual errors $\norm{u-u_X}{H^1(\Omega)}$
and $\norm{\lambda-\lambda_k}{L_2(\Gamma)}$ on a double logarithmic scale.
For the latter error, which is measured in the $L_2$ rather than $H^{-1/2}$-norm,
we expect a reduced convergence like $k^{1/2}$. Both expected error terms,
$h_X$ (labeled as $h$) and $k^{1/2}$, are also given in the plots
(multiplied by $10$ to shift them closer to the corresponding error curves).
Figures~\ref{fig_error_r02_k0} and \ref{fig_error_r02_k1} show the results
for $r=0.2$ with $\tau=1.5$ and $\tau=2.5$,
respectively. Figures~\ref{fig_error_r01_k0} and \ref{fig_error_r01_k1} show the
corresponding results for reduced radius $r=0.1$.
In all the cases there is some pre-asymptotic
range and the errors behave as expected for larger number of unknowns.

\begin{figure}[htb]
\begin{center}
\includegraphics[width=0.8\textwidth]{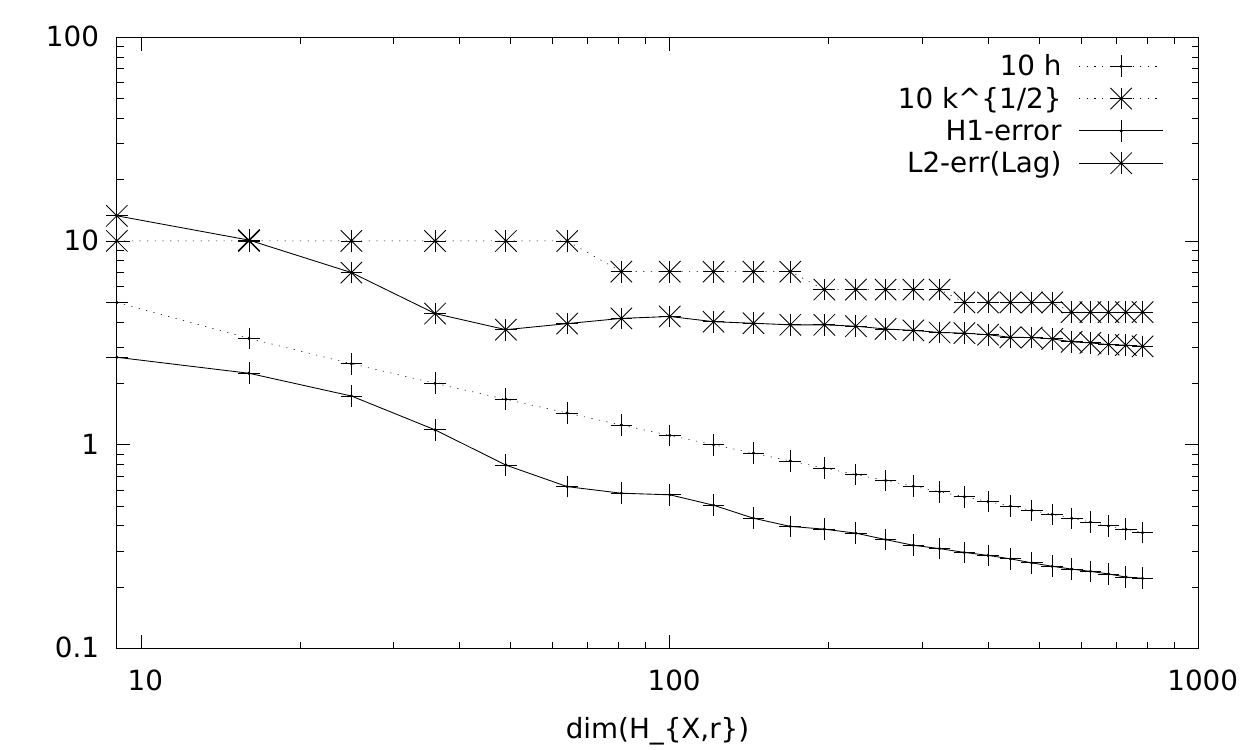}
\end{center}
\caption{Errors for $r=0.2$ and $\tau=1.5$.}
\label{fig_error_r02_k0}
\end{figure}

\begin{figure}[htb]
\begin{center}
\includegraphics[width=0.8\textwidth]{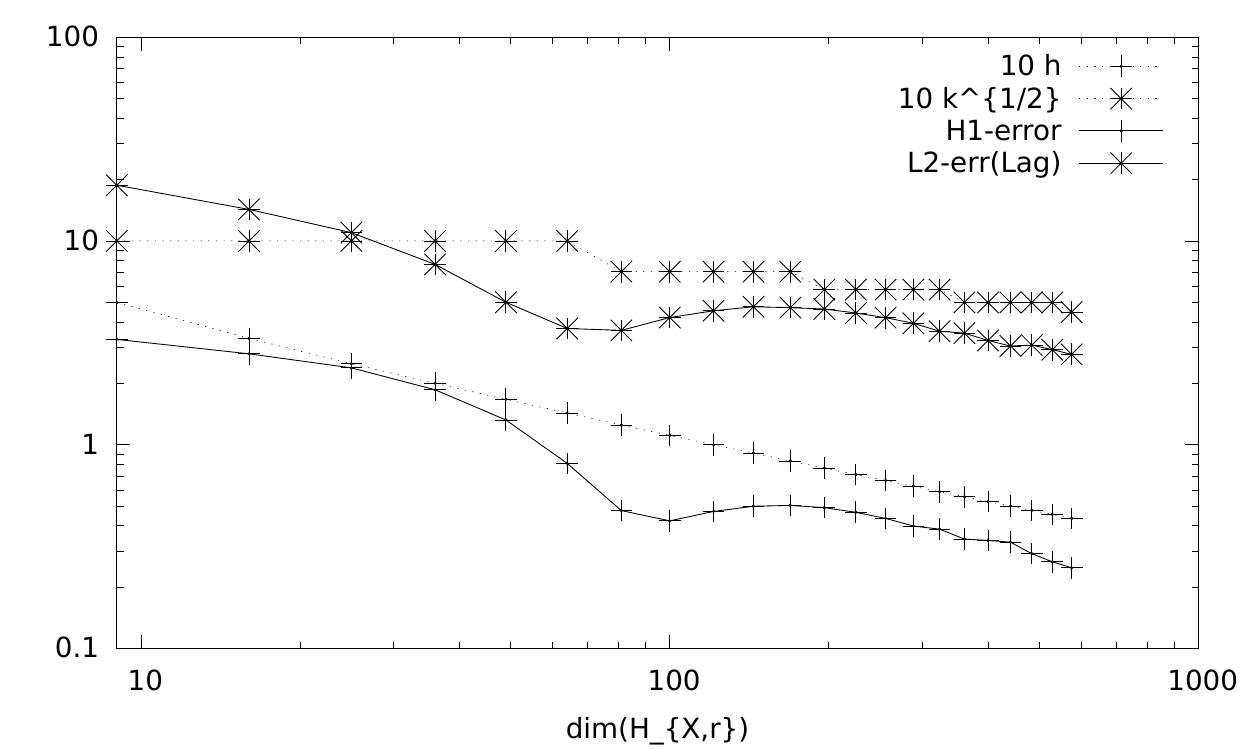}
\end{center}
\caption{Errors for $r=0.2$ and $\tau=2.5$.}
\label{fig_error_r02_k1}
\end{figure}

\begin{figure}[htb]
\begin{center}
\includegraphics[width=0.8\textwidth]{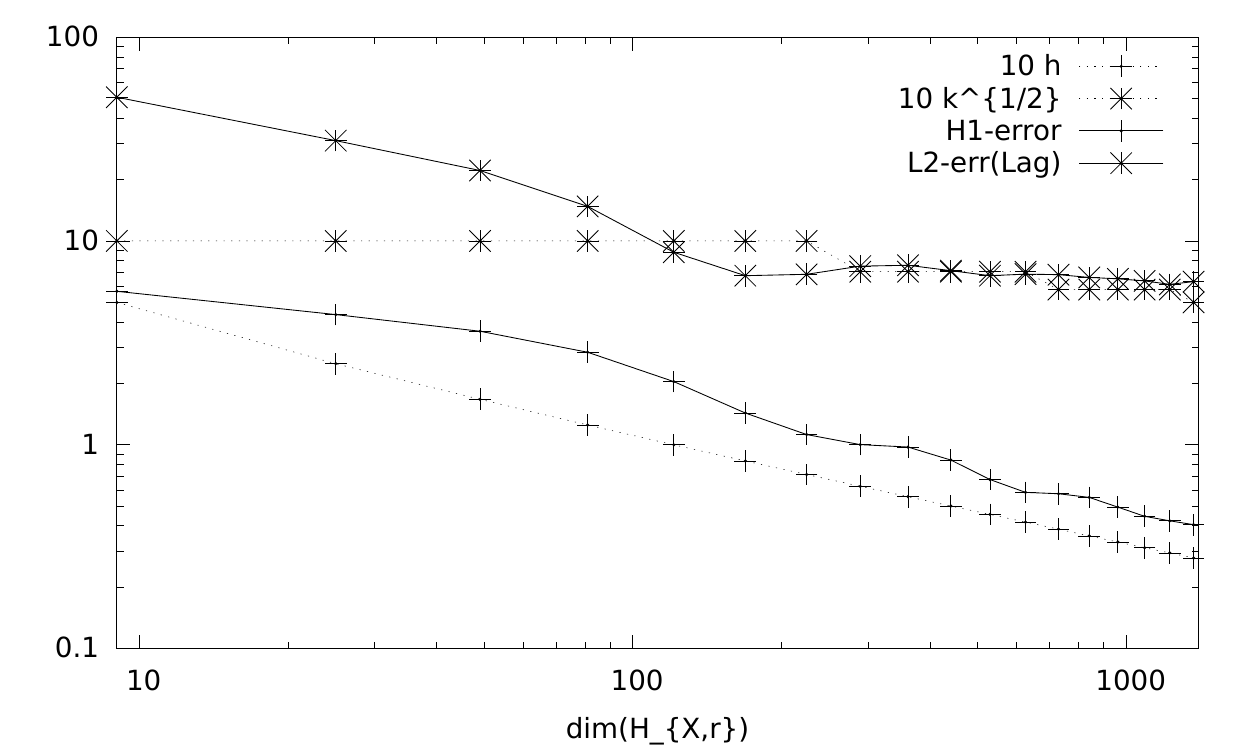}
\end{center}
\caption{Errors for $r=0.1$ and $\tau=1.5$.}
\label{fig_error_r01_k0}
\end{figure}

\begin{figure}[htb]
\begin{center}
\includegraphics[width=0.8\textwidth]{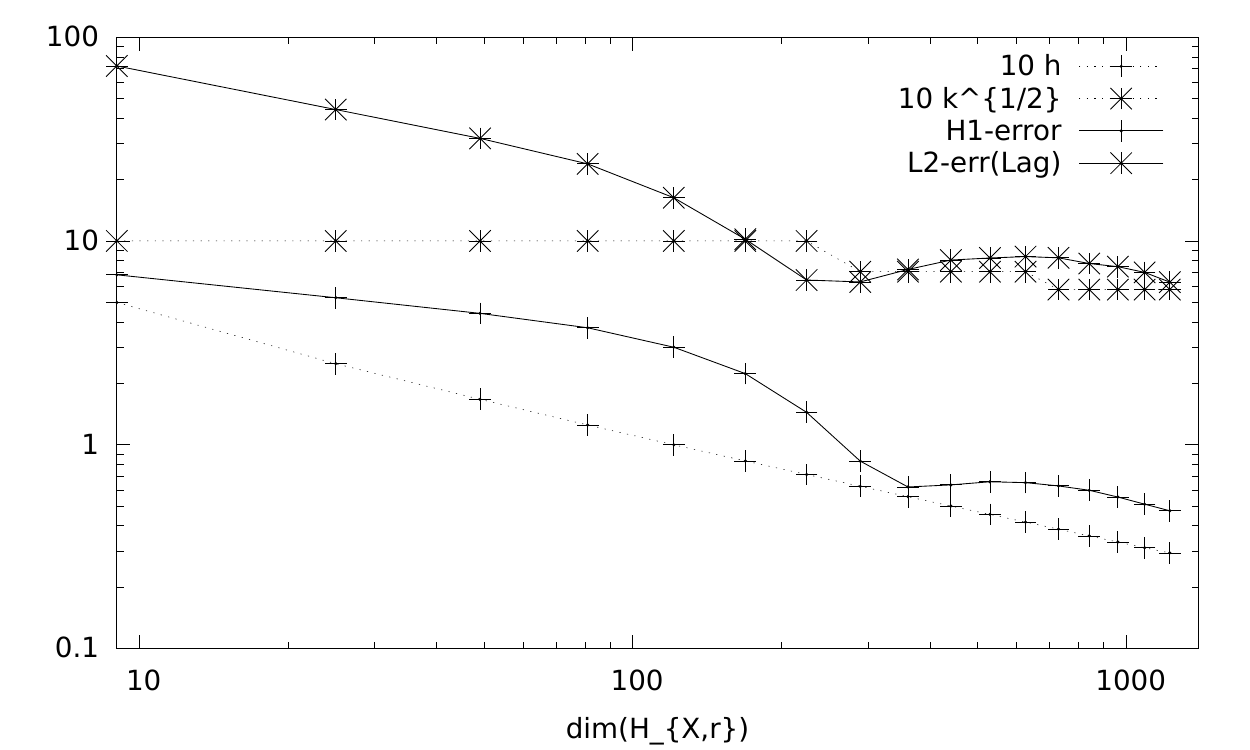}
\end{center}
\caption{Errors for $r=0.1$ and $\tau=2.5$.}
\label{fig_error_r01_k1}
\end{figure}

\def\cprime{$'$}

\end{document}